\documentclass[11pt]{amsart}
\usepackage{cprotect}

\usepackage[cmyk]{xcolor}
\definecolor{R}{cmyk}{0,1,1,0.0} 
\definecolor{Y}{cmyk}{0,0,.6,0} 
\definecolor{G}{cmyk}{.4,0,.8,0} 
\definecolor{B}{cmyk}{.8,0.4,0,0} 

\usepackage{tikz}
\usepackage{tikz-3dplot}
\usetikzlibrary{calc}
\tikzstyle{my state}=[rectangle,draw]
\tikzstyle{gray step}=[rectangle,draw]
\tikzstyle{Atrans}=[blue,ultra thick,->]
\tikzstyle{Ttrans}=[red,dashed,ultra thick,->]
\tikzstyle{box}=[rectangle,draw=black!50,fill=white,thick,inner sep=0pt,minimum size=4.5mm]

\usepackage{hyperref}
\usepackage{url}
\usepackage{mathrsfs}
\usepackage{breqn}
\usepackage[utf8]{inputenc}
\usepackage{amsmath}
\usepackage{amsthm}
\usepackage{amsfonts}

\usepackage{listings}
\definecolor{deepblue}{rgb}{0,0,0.5}
\definecolor{deepred}{rgb}{0.6,0,0}
\definecolor{deepgreen}{rgb}{0,0.5,0}
\DeclareFixedFont{\ttb}{T1}{txtt}{bx}{n}{12} 
\DeclareFixedFont{\ttm}{T1}{txtt}{m}{n}{12}  

\lstdefinestyle{mystyle}{
language=Python,
basicstyle=\ttm,
morekeywords={self},              
keywordstyle=\ttb\color{deepblue},
emph={MyClass,__init__},          
emphstyle=\ttb\color{deepred},    
stringstyle=\color{deepgreen},
frame=tb,                         
showstringspaces=false,
showspaces=false,
showtabs=true,
tabsize=2,
tab=\smash{\rule[-.2\baselineskip]{.4pt}{\baselineskip}\kern.5em}
}
\newtheorem{algorithm}{Algorithm}
\newtheorem*{theorem}{Theorem}
\newtheorem*{lemma}{Lemma}
\usepackage[capitalise]{cleveref}
\usepackage{graphicx}
\newcommand{\tmpbrk}[2]{}

\def\sign{\mathop{\mathrm{sign}}}
\author{Hans Jakob Rivertz}
\title[Multiset permutation generation]{Multiset permutation generation by transpositions}

\begin{document}

\begin{abstract}
This paper proposes a new algorithm for generating all permutations of multisets.
The method uses transpositions only and
adjacent transpositions are favoured.
The algorithm requires a strong homogeneous transposition condition:
non-adjacent transpositions are allowed only if all elements between the two permuted elements are equal to the smallest of those two elements.
The storage required by the algorithm is small. 
\end{abstract}
\maketitle
\section{Introduction}
\subsection{General introduction}
Combinatorics is considered as an independent discipline in mathematics and its history goes back several thousand years.
It is the origin of many of the branches of modern mathematics~\cite{Rota1969}.
The results and methods of combinatorics are also used in fields of mathematics with origin independent to combinatorics.
Problems of combinatoric nature exist in many of the several branches of mathematics today.
Differential geometry is not an exception.
It has even a sub field called combinatorial differential geometry.
This paper contains a solution to a combinatoric problem that arose in an ongoing research in differential geometry.
\cref{sec:motivation} gives a brief description of this problem.

Permutations and combinations are related through permutations of multisets.
In a set the elements are unique, whilst in a multiset, elements can occur several times.
A multiset can be described by a list of $r$ symbols $1,2,\ldots,r$, and their multiplicities $m_1,m_2,\ldots,m_r$.
The set $\{1,2,\ldots,r\}$ is a special case of a multiset.

A $k$-combination of elements in a set $S=\{1,2,3,\ldots,n\}$ is a subset $A$ of $S$ of size $k$.
Representing $A$ as a string of ones and zeros (a bit-string), where a one on the place $i$ says that $i$ is in $A$ and a zero on place $j$ tells that $j$ is not in $A$, gives a one-to-one correspondence between subsets of $S$ and bit-strings of length $n$.
Since $A$ has $k$ elements, the corresponding bit string has $k$ ones and $n-k$ zeros.
This is a permutation of the multiset consisting of $k$ ones and $n-k$ zeros.
Thus, generating all $k$-combinations is a special case of listing all permutations of the multiset $1^k0^{n-k}$.
The super script $k$ indicates there are $k$ copies of the element 1 in the multiset.

By carefully designing the $k$-combination generating algorithm, it can be used to generate all permutations of any multiset.
In fact, in a multiset consider the element 1 and treat the remaining elements as 0.
We require two conditions.
\begin{itemize}
\item[1.] \textbf{The transposition condition} requires that two successive permutations differs only in two bits.
\item[2.] \textbf{The zero stability condition} require that the order of the zeros is not changed.
\end{itemize}
For algorithms generating combinations, the first condition is equivalent to the strong minimal change property.
There is a rich literature considering those algorithms. A selection of algorithms can be found in~\cite{Korsh:2000,Kruchinin:2022,Takaoka1999_a,Takaoka1999,Takaoka2006,Torres2019}.

The second condition is similar to the homogeneous transpositions condition in~\cite{Ruskey2003}.
There are many algorithms~\cite{Knuth:2005:ACPb,Ruskey2003,Ruskey2005,Ruskey2009,Sawada2021,Williams2022,Williams2009} that meet other requirements.
These algorithms however, are not in the scope of this paper.

\subsection{Motivating introduction}
\label{sec:motivation}
The algorithm in the presented paper was developed to be used in making some huge polynomials.
The variables in these polynomials are the coefficients of the Riemann curvature tensor.
The polynomials are defined by using Young tableaux, which are diagrams filled with numbers.
The forms of the diagrams are defined by partitions.
A partition of $n$ is a sum $\lambda_1+\lambda_2+\cdots+\lambda_k=n$, where the terms are positive and non-increasing.
There are $\lambda_i$ boxes in row $i$ and the numbers $1,2,\ldots,n$ are distributed in the boxes. The following Young tableau $B_{2,2,2,2}$ is defined from the partition $2+2+2+2=8$.
$$
B_{2,2,2,2}=
\begin{tikzpicture}[
baseline=($(A.base)!.5!(H.base)$),
auto,node distance=4.8mm,semithick,inner sep=2pt,bend angle=45]
\node [box] (A) {1};
\node [box] (B) [below of=A] {2};
\node [box] (C) [right of=A] {3};
\node [box] (D) [below of=C] {4};
\node [box] (E) [below of=B] {5};
\node [box] (F) [below of=E] {6};
\node [box] (G) [below of=D] {7};
\node [box] (H) [below of=G] {8};
\end{tikzpicture}
$$
This Young Tableau was used by Agaoka~\cite{agaoka85} to find a necessary condition on isometric immersions of four-dimensional Riemann spaces into the Euclidean space of dimension five.

A Young tableau has a corresponding Young symmetrizer which is used to produce symmetric polynomials.
The symmetrizer is defined through subgroups of the symmetric group $\mathfrak{S}_{8}$.
The vertical subgroup $\mathfrak{V}_{B}$ consists of all elements in  $\mathfrak{S}_{8}$ that preserves the columns of the Young Tableau $B$.
It's size is $(4!)^2=576$.
The horizontal subgroup $\mathfrak{H}_{B}$ consists of all elements in  $\mathfrak{S}_{8}$ that preserves the rows of the Young Tableau $B$.
The size of $\mathfrak{H}_{B}$ is $(2!)^4=16$.
The polynomial is defined as
\begin{dmath*}
\sum_{\sigma^{-1}\in\mathfrak{V}_{B}}
\sign(\sigma)
R(Y_{\sigma(1)},
  Y_{\sigma(2)},
  Y_{\sigma(3)},
  Y_{\sigma(4)})\tmpbrk{\cdot}{}
R(Y_{\sigma(5)},
  Y_{\sigma(6)},
  Y_{\sigma(7)},
  Y_{\sigma(8)}),
\end{dmath*}
where $Y_i=X_{\mathscr{R}(i)}$. $\mathscr{R}(i)$ is the row where the number $i$ is in the Young tableau $B$. The symbols $$X_1,X_2,X_3,X_4$$ denotes a basis for $\mathbb{R}^4$.
The multivariate function $R$ is a $(0,4)$-tensor, that is a  real valued multilinear function with arguments in $\mathbb{R}^4$.
Let $R$ have the symmetries
\begin{itemize}
\item[$S_1$.] $R(X,Y,Z,W)=-R(Y,X,Z,W)$,
\item[$S_2$.] $R(X,Y,Z,W)=-R(X,Y,W,Z)$,
\item[$S_3$.] $R(X,Y,Z,W)=R(Z,W,X,Y)$, and
\item[$S_4$.] $R(X,Y,Z,W)+R(Z,X,Y,W)+R(Y,Z,X,W)=0$. (The first Bianchi identity.)
\end{itemize}
$(0,4)$-tensors $R$ that have these symmetries are called curvature like tensors since they share the same symmetries as the Riemann curvature tensor.
From the first symmetry, we have for example
$$R(Y_{\sigma(1)},
  Y_{\sigma(2)},
  Y_{\sigma(3)},
  Y_{\sigma(4)})
=-
R(Y_{\sigma(2)},
  Y_{\sigma(1)},
  Y_{\sigma(3)},
  Y_{\sigma(4)}).$$
Therefore, the two permutations which only differs by a transposition of 1 and 2, gives the same term in the polynomial.
The same is true for the pairs $(3,4)$, $(5,6)$, and $(7,8)$.
Instead of summing over $\mathfrak{V}_B$, we will sum over all permutations of the multiset obtained by identifying elements in the same pair.
The number of terms in the sum has been reduced by a factor of $2^4$ to $\frac{(4!)^2}{16}=36$.
In addition to the four symmetries of $R$, there is a symmetry in the polynomial:
interchanging the factors in each term is in fact equivalent to a permutation in $\mathfrak{V}_B$.
This permutation is $(1,5)(2,6)(3,7)(4,8)$.
This paper considers first the symmetry types $S_1$ and $S_2$ mainly.
The other symmetries are applied in \cref{sec:gen_pol}.
 
As can be seen, the symmetries in the tensors reduce the complexity of the calculations.
Using transpositions in the generation of the permutations also makes it trivial to calculate the factor $\sign(\sigma)$; it changes sign for each step.
\section{Related algorithms and theory}
\subsection{The Steinhaus-Johnson-Trotter algorithm}
The Steinhaus-Johnson-Trotter algorithm~\cite{Selmer_M_Johnson:1963,Trotter:1962} generates all permutations of $n$ distinct elements by adjacent transpositions only.
As Donald Knuth points out on page 321, volume 4 in The `Art of Computer Programming'~\cite{Knuth:2005:ACPb}, it is impossible to generate all permutations even, the simple multiset \{1,1,2,2\}.
The proof is an application of graph theory.
The edges in the following diagram are adjacent transpositions.
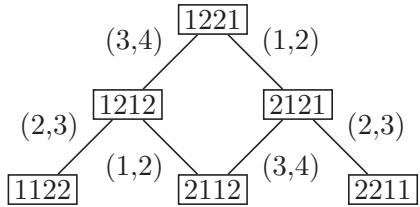
\begin{center}
\begin{tikzpicture}[auto,node distance=1.6cm,semithick,inner sep=2pt,bend angle=45]
\node [my state] (A) {1122};
\node [my state] (B) [above right of=A] {1212};
\node [my state] (C) [above right of=B] {1221};
\node [my state] (D) [below right of=B] {2112};
\node [my state] (E) [below right of=C] {2121};
\node [my state] (F) [below right of=E] {2211};
\path (A) edge node {(2,3)} (B)
      (B) edge node {(3,4)} (C)
      (C) edge node {(1,2)} (E)
      (D) edge node {(1,2)} (B)
      (E) edge node {(3,4)} (D)
      (E) edge node {(2,3)} (F);
\end{tikzpicture}
\end{center}
There is no Hamilton path in the graph above, so there is no way of generating all permutations of $\{1,1,2,2\}$ by successive adjacent transpositions.
\subsection{Gray codes}
For combinations, a Gray code is a complete list of all $k$-combinations of $n$ elements listed by successive transpositions of the bit-string representation.
On page 127 in~\cite{Ruskey2003}, Ruskey presents the following algorithm.
The list $C(n,k)$ of all permutations of $1^k0^{n-k}$ is given recursively by
$$C(n,k)=\begin{cases}
0^n&\mbox{if $k=0$}\\
1^n&\mbox{if $k=n$}\\
\begin{bmatrix}
C(n-1,k)\cdot 0\\
\overline{C(n-1,k-1)}\cdot 1
\end{bmatrix}
&\mbox{if $0<k<n$}
\end{cases}.$$
The line over $C(n-1,k-1)$ indicates reversion of the list and the dot in $L\cdot c$ means concatenating the character $c$ to all of the elements in the list $L$.
The production of $C(4,2)$ is described in the following diagram.
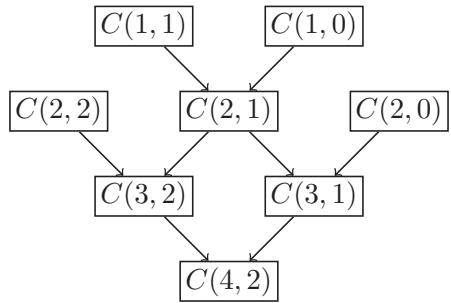
\begin{center}
\begin{tikzpicture}[auto,node distance=1.6cm,semithick,inner sep=2pt,bend angle=45]
\node[gray step] (A) {$C(4,2)$};
\node[gray step] (B) [above left  of=A] {$C(3,2)$};
\node[gray step] (C) [above right of=A] {$C(3,1)$};
\node[gray step] (D) [above left  of=B] {$C(2,2)$};
\node[gray step] (E) [above right of=B] {$C(2,1)$};
\node[gray step] (F) [above right of=C] {$C(2,0)$};
\node[gray step] (G) [above left of=E] {$C(1,1)$};
\node[gray step] (H) [above right of=E] {$C(1,0)$};
\path (H) edge[->] (E)
      (G) edge[->] (E)
      (F) edge[->] (C)
      (E) edge[->] (C)
      (E) edge[->] (B)
      (D) edge[->] (B)
      (C) edge[->] (A)
      (B) edge[->] (A);
\end{tikzpicture}
\end{center}
This gives the listing $1100,0110,1010,0011,0101,1001$.
In the adjacent transposition diagram, this list is represented as following.
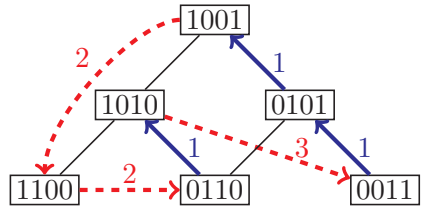
\begin{center}
\begin{tikzpicture}[auto,node distance=1.6cm,semithick,inner sep=2pt,bend angle=45]
\node [my state] (A) {1100};
\node [my state] (B) [above right of=A] {1010};
\node [my state] (C) [above right of=B] {1001};
\node [my state] (D) [below right of=B] {0110};
\node [my state] (E) [below right of=C] {0101};
\node [my state] (F) [below right of=E] {0011};
\path (A) edge[Ttrans] node[above] {2} (D)
      (D) edge[Atrans] node[right=1mm] {1} (B)
      (B) edge[Ttrans] node[right=4mm] {3} (F)
      (F) edge[Atrans] node[right=1mm] {1} (E)
      (E) edge[Atrans] node[right=1mm] {1} (C)
      (A) edge (B)
      (B) edge (C)
      (E) edge (D);
\draw [Ttrans](C) .. controls +(180:1) and +(90:1) .. node[above left] {2}  (A);
\end{tikzpicture}
\end{center}
The thick arrows correspond to adjacent transpositions.
The non-adjacent transpositions are shown as dashed arrows.
The numbers on the arrows show the distance between the changing bits.
Notice that there is a transposition from the last bit string to the first bit string.
In fact, the algorithm produces a loop of transpositions, where the list can be repeated.
In the process of the above example, if the zeros and the ones where marked, their order will change.
However, when flipping back to start, the order is unchanged.
For example the marked listing $C'(4,2)$ is
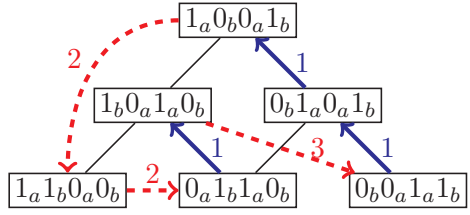
\begin{center}
\begin{tikzpicture}[auto,node distance=1.6cm,semithick,inner sep=2pt,bend angle=45]
\node [my state] (A1) {$1_a1_b0_a0_b$};
\node [my state] (B1) [above right of=A1] {$1_b0_a1_a0_b$};
\node [my state] (C1) [above right of=B1] {$1_a0_b0_a1_b$};
\node [my state] (A2) [below right of=B1] {$0_a1_b1_a0_b$};
\node [my state] (B2) [below right of=C1] {$0_b1_a0_a1_b$};
\node [my state] (A3) [below right of=B2] {$0_b0_a1_a1_b$};
\path (A1) edge[Ttrans] node[above] {2} (A2)
      (A3) edge[Atrans] node[right=1mm] {1} (B2)
      (A2) edge[Atrans] node[right=1mm] {1} (B1)
      (B1) edge[Ttrans] node[right=3mm] {3} (A3)
      (B2) edge[Atrans] node[right=1mm] {1} (C1)
      (A1) edge (B1)
      (A2) edge (B2)
      (B1) edge (C1);
\draw [Ttrans](C1) .. controls +(180:2) and +(90:1) .. node[above left] {2} (A1);
\end{tikzpicture}
\end{center}
In general, the order of the marked ones will change as is seen in the following example:
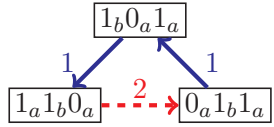
\begin{center}
\begin{tikzpicture}[auto,node distance=1.6cm,semithick,inner sep=2pt,bend angle=45]
\node [my state] (AL) {$1_a1_b0_a$};
\node [my state] (B) [above right of=AL] {$1_b0_a1_a$};
\node [my state] (AR) [below right of=B] {$0_a1_b1_a$};
\path (AL) edge[Ttrans] node[above] {2} (AR)
      (AR) edge[Atrans] node[right=2mm] {1} (B)
      (B) edge[Atrans] node[left=2mm] {1} (AL);
\end{tikzpicture}
\end{center}
However, it is known that the order of the marked zeros will not change in the closed loop of the algorithm.
Define the marked zeros version of the algorithm recursively by
$$C'(n,k)=\begin{cases}
0_00_1\cdots0_{n-1}&\mbox{if $k=0$}\\
1^n&\mbox{if $k=n$}\\
\begin{bmatrix}
C'(n-1,k)\cdot 0_{n-k-1}\\
\overline{C'(n-1,k-1)}\cdot 1
\end{bmatrix}
&\mbox{if $0<k<n$}
\end{cases}.$$
The line over $C'(n-1,k-1)$ indicates reversion of the list and subtracting one from each index modulo $n-k$.
\begin{lemma}
For all $n>0$ and for all $k$ so that $0\leq k\leq n$ the following is true.
\begin{itemize}
\item[1.] The first element of the list $C'(n,k)$ is
$$1^k0_00_1\cdots0_{n-k-1}.$$
\item[2.] When $0<k<n$, the last element of the list $C'(n,k)$ is
$$1^{k-1}0_{n-k-1}0_00_1\cdots0_{n-k-2}1.$$
\item[3.] Each two successive elements in the list $C'(n,k)$ differ in exactly two symbols, a marked zero and a one.
\end{itemize}
\end{lemma}
Notice that the first and last elements in the list differs in only in the $k$-th and $n$-th bit when the list has more than one element.
\begin{proof}[Proof of the lemma]
The proof is by induction on $n$.
Notice that the lemma is true if $n=1$, $k=0$, or $k=n$.
Assume the lemma is true for $n<r$, for some $r>1$. It is sufficient to prove the lemma for $n=r$ and $0<k<n$.
First, from the definition one can see that the first element in $C'(r,k)$ is equal to the first element in $C'(r-1,k)\cdot 0_{r-k-1}$.
That is
$$1^k0_00_1\cdots0_{n-k-2}\cdots0_{n-k-1}.$$
Second, from the definition, the last element in $C'(r,k)$ is the first element of $C'(r-1,k-1)\cdot 1$ with one subtracted from the indices modulo $r-k$.
That is
$$1^{k-1}\cdots0_{r-k-1}0_00_1\cdots0_{r-k-2}1.$$
Finally, it is sufficient to prove the third point in the lemma for the last element in $C'(r-1,k)\cdot 0_{r-k-1}$ and the first element in $\overline{C'(r-1,k-1)}\cdot 1$.
The last element in $C'(r-1,k)\cdot 0_{r-k-1}$ is
$$1^{k-1}0_{r-k-1}0_00_1\cdots0_{r-k-2}1.$$
The first element in $\overline{C'(r-1,k-1)}\cdot 1$ is the last element in $C'(r-1,k-1)\cdot 1$ with the one subtracted indices modulo $r-k$.
That is
$$1^{k-2}0_{r-k-2}0_{r-k-1}0_00_1\cdots0_{r-k-3}11.$$
Notice that these two displayed elements differs in 
place $k-1$ and $r-1$ only.
The lemma follows from induction.
\end{proof}
An implementation of the definition of $C(n,k)$ requires memory allocated for the entire list.
Two alternative algorithms are given in chapter 5.3 on the pages 128 and 129 in~\cite{Ruskey2003}.
In the algorithm on page 128, each element is obtained from the previous by a transposition.
The bit-string representation of a $k$-combination is equivalent to an ascending ordered sequence $p_1p_2\cdots p_k$ of the position of the ones in the bit-string.
The second algorithm gives the list $C(n,k)$ in terms of position sequences.
Mütze has written a recent survey~\cite{Mutze:2022} on combinatorial Gray codes.

\subsection{The piano chords algorithm}
Eades and McKay~\cite{Eades:1984} found an algorithm that theoretically goes through all chords using $k$ fingers on a piano with $n$ keys without crossing any finger and moving only one finger at the time.
It is theoretical because no human has so large or flexible hands that they can play the chords in this list.
The algorithm is described as following in~\cite{Ruskey2003}.
$$E(n,k)=\begin{cases}
0^n&\mbox{ if $k=0$}\\
[10^{n-1},010^{n-2},\ldots,0^{n-1}1]^T&\mbox{ if $k=1$}\\
\begin{bmatrix}
E(n-1,k)\cdot0\\
\overline{E(n-2,k-1)}\cdot01\\
E(n-2,k-2)\cdot11
\end{bmatrix}&\mbox{ if $1<k<n$}\\
1^n&\mbox{ if $k=n$}
\end{cases}
$$
For $n=6$ and $k=3$, this gives the list given in \cref{tab:0}. There are 16 adjacent transpositions in the generation, two transpositions of distance two and one transposition of length three. The total distance moved in this algorithm is 23. 
\begin{table}[htb!]
\begin{center}
\begin{tabular}{|rlr|rlr|}\hline
&state&trans.&&state&trans.\\\hline
 1&\verb|111000|&     &11&\verb|001101|&(5,6)\\
 2&\verb|110100|&(3,4)&12&\verb|100101|&{\bf(1,3)}\\
 3&\verb|101100|&(2,3)&13&\verb|010101|&(1,2)\\
 4&\verb|011100|&(2,3)&14&\verb|011001|&(3,4)\\
 5&\verb|011010|&(4,5)&15&\verb|101001|&(1,2)\\
 6&\verb|101010|&(1,2)&16&\verb|110001|&(2,3)\\
 7&\verb|110010|&(2,3)&17&\verb|100011|&{\bf(2,5)}\\
 8&\verb|100110|&{\bf(2,4)}&18&\verb|010011|&(1,2)\\
 9&\verb|010110|&(1,2)&19&\verb|001011|&(2,3)\\
10&\verb|001110|&(2,3)&20&\verb|000111|&(3,4)\\\hline
\end{tabular}
\end{center}
\caption{The piano chord algorithm by Eades and McKay applied with $n=6$ and $k=3$. The third and sixth columns show the transposition that created the state from the preceding state.}
\label{tab:0}
\end{table}
\section{The new algorithm}
\subsection{Definition and proofs}
The algorithm start with $k$ identical elements marked with the letter $x$ placed in location $1,2,\ldots,k$ along a linear array with $n$ places
$$
\underbrace{x\cdots x}_{k\mbox{ copies}}\!\!o\cdots o.
$$
The empty places are marked with the letter $o$.
The algorithm describes a process of transpositions that ends up in the state
$$
o\cdots o\!\!\underbrace{x\cdots x}_{k\mbox{ copies}},
$$
and in the process going through all the $\binom nk$ configurations of $k$ number of $x$'s and $n-k$ number of $o$'s exactly once.
\begin{algorithm}
\label{alg:1}
Given an array of $k$ number of $x$'s distributed in an array with $n$ positions.
The $x$'s are given an orientation coded by left ($<$) or right ($>$).
The empty places are represented by $o$ and they have no orientation.
The $x$'s are therefore redundant and are represented by the symbol assigning its orientation.
The algorithm moves the elements and changes their orientation by the following set of rules.
\begin{itemize}
\item[1.] Output the state. The active element is set to the element furthest to the right.
\item[2.] If the active element faces an adjacent empty place, it interchanges with that empty place.
$$\ldots >o\ldots\mapsto\ldots o>\ldots\mbox{ or }
\ldots o<\ldots\mapsto\ldots <o\ldots$$
\noindent Go back to step 1.
\item[3.] If the active element faces left and faces an empty place to its left and all places between are occupied by elements facing to the left, it permutes with that empty place. The elements between change orientation to right ($>$). Go back to the first step.
\item[4.] If the active element faces an empty place to its right, it will permute with the closest empty place to its right if and all elements between are oriented to the right ($>$). The algorithm goes to step 1.
\item[5.] The active element changes its orientation. The nearest element to the left of the active element becomes the active element and the algorithm goes to step 2.
\item[6.] The algorithm terminates if there are no elements to the left of the active element.
\end{itemize}
The algorithm runs through all combinations with the following initial states.
\begin{itemize}
\item[a.]If the initial state is
$$
\underbrace{>\cdots >}_{k\mbox{ copies}}\!o\cdots o,
$$
the procedure defined in point 1 to 6 describes an algorithm that will terminate with the last output on the form
$$
o\cdots o>\!\!\!\!\!\!\underbrace{*\cdots *}_{k-1\mbox{ copies}}\!\!\!\!,
$$
where $*$ can be any of the two symbols $<$ and $>$.
\item[b.]If the initial state is the negative of the last output in case a:
$$
o\cdots o<\!\!\!\!\!\!\underbrace{*\cdots *}_{k-1\mbox{ copies}}\!\!\!\!,
$$
the procedure defined in point 1 to 6 describes an algorithm that will terminate with the last output on the form
$$
\underbrace{<\cdots <}_{k\mbox{ copies}}\!o\cdots o.
$$
\end{itemize}
\end{algorithm}
A proof of the Algorithm is given at the end of this section.
\cref{tab:1} shows the process with three elements and six places.
The first column counts the set permutations,
the second column shows the state, and the third column indicates the rules that produced that state.
Notice that \cref{tab:1} shows all $\binom 63=20$ combinations of three unmarked elements in six places. The orientations are for book-keeping only.
\begin{table}[htb!]
\begin{center}
\begin{tabular}{|rlr|rlr|}\hline
&state&rules&&state&rules\\\hline
1&\verb|>>>ooo|&&11&\verb|o>oo<>|&6,6,2\\
2&\verb|>>o>oo|&2&12&\verb|o>o<o<|&6,2\\
3&\verb|>>oo>o|&2&13&\verb|o>o<<o|&2\\
4&\verb|>>ooo>|&2&14&\verb|o><>oo|&\bf3\\
5&\verb|>o>oo<|&6,2&15&\verb|o><o>o|&2\\
6&\verb|>o>o<o|&2&16&\verb|o><oo>|&2\\
7&\verb|>o><oo|&2&17&\verb|oo>>o<|&6,6,\bf4\\
8&\verb|>oo>>o|&6,\bf4&18&\verb|oo>><o|&2\\
9&\verb|>oo>o>|&2&19&\verb|oo>o>>|&6,\bf4\\
10&\verb|>ooo><|&6,2&20&\verb|ooo><<|&6,6,2\\\hline
\end{tabular}
\end{center}
\caption{The table shows the algorithm applied on the state \texttt{\symbol{62}\symbol{62}\symbol{62}ooo}.}
\label{tab:1}
\end{table}

\cref{tab:2} shows the algorithm when the start configuration is \verb|ooo<>>|.
\begin{table}[htb!]
\begin{center}
\begin{tabular}{|rlr|rlr|}\hline
&state&rules&&state&rules\\\hline
1&\verb|ooo<>>|&&11&\verb|<ooo<>|&6,6,2\\
2&\verb|oo<o<<|&6,6,2&12&\verb|<oo<o<|&6,2\\
3&\verb|oo<<>o|&\bf3&13&\verb|<oo<<o|&2\\
4&\verb|oo<<o>|&2&14&\verb|<o<>oo|&\bf3\\
5&\verb|o<>oo<|&6,\bf3&15&\verb|<o<o>o|&2\\
6&\verb|o<>o<o|&2&16&\verb|<o<oo>|&2\\
7&\verb|o<><oo|&2&17&\verb|<<ooo<|&6,2\\
8&\verb|o<o>>o|&6,\bf4&18&\verb|<<oo<o|&2\\
9&\verb|o<o>o>|&2&19&\verb|<<o<oo|&2\\
10&\verb|o<oo><|&6,2&20&\verb|<<<ooo|&2\\\hline
\end{tabular}
\end{center}
\caption{The table shows the algorithm applied on the state \texttt{ooo\symbol{60}\symbol{62}\symbol{62}}, the negative state of the end result in \cref{tab:1}.}
\label{tab:2}
\end{table}
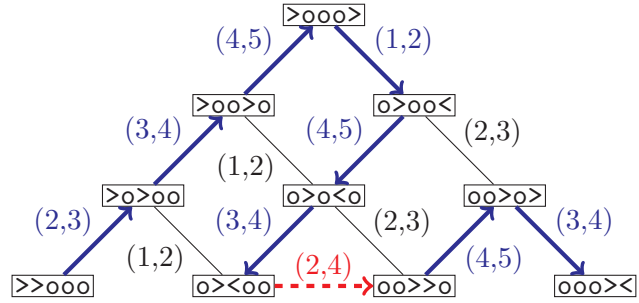
\begin{figure}[htb!]
\begin{center}
\begin{tikzpicture}[auto,node distance=1.7cm,thin,inner sep=1pt,bend angle=45]
\node [my state] (A) {\verb|>>ooo|};
\node [my state] (B) [above right of=A] {\verb|>o>oo|};
\node [my state] (C) [above right of=B] {\verb|>oo>o|};
\node [my state] (D) [below right of=B] {\verb|o><oo|};
\node [my state] (E) [below right of=C] {\verb|o>o<o|};
\node [my state] (F) [above right of=C] {\verb|>ooo>|};
\node [my state] (G) [above right of=E] {\verb|o>oo<|};
\node [my state] (H) [below right of=E] {\verb|oo>>o|};
\node [my state] (I) [above right of=H] {\verb|oo>o>|};
\node [my state] (J) [below right of=I] {\verb|ooo><|};
\path (A) edge[Atrans] node {(2,3)} (B)
      (B) edge[Atrans] node {(3,4)} (C)
      (C) edge[Atrans] node {(4,5)} (F)
      (F) edge[Atrans] node {(1,2)} (G)
      (G) edge[Atrans] node[above left] {(4,5)} (E)
      (E) edge[Atrans] node[above left] {(3,4)} (D)
      (D) edge[Ttrans] node {(2,4)} (H)
      (H) edge[Atrans] node[below right] {(4,5)} (I)
      (I) edge[Atrans] node {(3,4)} (J)
      (E) edge node {(2,3)} (H)
      (G) edge node {(2,3)} (I)
      (E) edge node {(1,2)} (C)
      (D) edge node {(1,2)} (B);
\end{tikzpicture}
\end{center}
\caption{The adjacent transpositions graph for permutations of the multiset \texttt{xxooo}. The bold arrows correspond to the adjacent transpositions in the algorithm. The dashed arrow shows the only case where the transposition fails to be adjacent.}
\end{figure}
Notice the symmetry of the algorithm in \cref{tab:1} and \cref{tab:2}.
This symmetry is independent of the number of places and arrows.
\begin{theorem}
For each string $S$ of the characters `\verb|o|', `\verb|<|', and `\verb|>|', define its negative $N(S)$ by keeping `\verb|o|', replacing `\verb|<|' with `\verb|>|', and replacing `\verb|>|' with `\verb|<|'.
For example $N(\texttt{o<oo<>})=\texttt{o>oo><}$.
Let $M(S)$ be the result of running one iteration of the algorithm.
Then for all strings $S$ formed by the characters `\verb.o.', `\verb.<.',and `\verb.>.' the following holds:
$$M\circ N\circ M(S)=N(S).$$
\end{theorem}
Consider an integer represented in the balanced ternary number format representation.
This is $$(t_nt_{n-1}\ldots t_2t_1t_0)_{\bar 3}=\sum_{i=0}^nt_i3^i,\,t_i\in\{-1,0,1\}, i=0,1,\ldots,n.$$
For instance, $42=(1\bar1\bar1\bar10)_{\bar3}$, where $\bar1$ is the digit $-1$.
If \verb|o| represents the digit zero, \verb|>| represents the digit one, and \verb|<| represents digit negative one, then $N(S)$ represents the negative value of $S$.
Thus, the notation $N(S)$ is chosen.
\begin{proof}[Proof of the theorem]
An unsuccessful iteration runs through all elements and changes their orientation. That is the same as the taking the negative.
Then $N\circ M(S)=S$.
A successful iteration interchanges an element with an empty space.
The proof is divided into three cases:
\begin{itemize}
\item If the interchange happens in step 2, all elements to the right of the current active element has changed orientation before the current element became active. These elements are changed back again with negation, and therefore are unchanged by $N\circ M$. The current active element is contained in a substring \texttt{>o} or \texttt{o<}. This is changed by $N\circ M$ to \texttt{o<} or \texttt{>o} respectively. No elements to the left of the active element is changed by $M$. Therefore, $N\circ M\circ N\circ M(S)=S$.
\item If the interchange happens in step 3, just as for step 2, all elements to the right of the active element changes orientation by $M$ and therefore unchanged by $N\circ M$. The active element is the element furthest to the right in a series $\texttt{o<}\ldots\texttt{<}$. This is changed to $\texttt{<>}\ldots\texttt{>o}$ by $M$ and then to $\texttt{><}\ldots\texttt{<o}$ by $N$.
The element furthest to the right of these elements will be the active element in the second application of $M$, but then step 4 in the algorithm applies. The string will change to $\texttt{o>}\ldots\texttt{>}$ and then to $\texttt{o<}\ldots\texttt{<}$ by $N$.
The elements to the left and to the right will not be changed by $N\circ M\circ N\circ M$, therefore $N\circ M\circ N\circ M(S)=S$.
\item If the interchange happens in step 4, all elements to the right of the active element has changed orientation. The active element is the left most element in a string $\texttt{>}\ldots\texttt{>o}$ before the action in step 4.
The action in step 4 changes this string to $\texttt{0>}\ldots\texttt{>}$.
Then $N$ changes this substring to $\texttt{0<}\ldots\texttt{<}$. The next time $M$ acts on this it will change the string to $\texttt{<>}\ldots\texttt{>0}$ by step 3. Then $N$ will change the string to $\texttt{><}\ldots\texttt{<0}$ this was the original state. The elements to the left and right are not changed by $N\circ M\circ N\circ M$, so $N\circ M\circ N\circ M(S)=S$.
\end{itemize}
This completes the proof.
\end{proof}
Consider the \textbf{reduction map} $r:S\mapsto r(S)$ that removes the rightmost character \verb|>| or \verb|<| from the string $S$.
For example if $S=\texttt{>oo>>o}$ then $r(S)=\texttt{>oo>o}$.
This map transforms \cref{tab:1} into the shorter \cref{tab:3}.
Multiple recurring states are removed.
The numbering of the states in \cref{tab:1} is kept in \cref{tab:3}.
Notice that the same table would be obtained if the algorithm was used directly on the start configuration \verb|>>ooo|.
This observation is basis for the proof of the algorithm.

\begin{table}[htb!]
\begin{center}
\begin{tabular}{|rlr|rlr|}\hline
&state&rules&&state&rules\\\hline
1&\verb|>>ooo|&&12&\verb|o>o<o|&2\\
5&\verb|>o>oo|&2&14&\verb|o><oo|&2\\
8&\verb|>oo>o|&2&17&\verb|oo>>o|&6,4\\
10&\verb|>ooo>|&2&19&\verb|oo>o>|&2\\
11&\verb|o>oo<|&6,2&20&\verb|ooo><|&6,2\\\hline
\end{tabular}
\end{center}
\caption{This table shows the algorithm applied on the state \texttt{\symbol{62}\symbol{62}ooo}. It also shows the reduction of \cref{tab:1}.}
\label{tab:3}
\end{table}

\begin{proof}[Proof of \cref{alg:1}]
Consider any state $S'$ in the algorithm with $k=l$ and $n=m$.
This state can be extended to a state $S$ by adding an extra character \verb|>| or \verb|<|.
\begin{itemize}
\item[(R)] The character \verb|>| is placed just to the right of the element element furthest to the right different from \verb|o| in $S'$. After successive iterations with rule 2 we run through all states that reduces to $S'$.
In the next iteration, the element element furthest to the right in $S$ changes orientation and the element element furthest to the right in $S'$ becomes the new active element.
If $S=$ \verb|...>oooo>|, $S=$ \verb|...<o...o>| or $S=$ \verb|...<>| or $S=$ \verb|...>>|, the result of the iteration on $S$ will commute with the reduction. $M\circ r=r\circ M$.

\item[(L)] The character \verb|<| is appended at the right end of $S'$. That is $S=$ \verb|..>o..o<|.
Repeated iterations with step 2, goes through all states of $S$ with $r(S)=S'$ until we get one of the following states.
$S=$ \verb|..?><o..o| or $S=$ \verb|..<<o..o|.
The next iteration for the first case invokes step 6 once and step 4 once, the result is $M(S)=$ \verb|..?o>>o..o|.
Again $r\circ M=M\circ r$. In the second case step 3 or 6 is invoked and again the iteration $M$ commutes with reduction $r$.
\end{itemize}

The statement in the algorithm is trivially true for all $n$ if $k=1$.
Assume that the algorithm part a is true for $n=m$ and for $k=l$.
We will prove that the algorithm is true for $n=m+1$ and $k=l+1$.

From the start configuration $S_0=$ \verb|>...>>o...o| with $k=l+1$ and $n=m+1$, the element furthest to the left is moving to the right for each iteration in the algorithm until it comes to the right most place. Rule 2 in the algorithm is applied each time.
After $m-l$ iterations, the state is $S_{m-l}=$ \verb|>...>o...o>|.
The reduction map maps each of these states to the start configuration \verb|>...>o...o| for $k=l$ and $n=m$.
In the next iteration rule 6 is applied, the orientation of the active element changes to \verb|<| and the next element is the active element.
The state is now \verb|>...>o...o<| and rule 2 is invoked on the element at position $l$. 
That results in the state $S_{m-l+1}=$ \verb|>...>o>o...o<|.
The reduced map sends this state onto $S'_1=r(S_{m-l+1})= \verb|>...>o>o...o|$. This is equivalent with running the first iteration with $k=l$ and $n=m$.

In general consider any output state in the algorithm on one of the two forms:
\begin{itemize}
\item[(F)] $S=$\verb|?...*>o...o| with possible zero \verb|o|'s at the right end.
\item[(B)] $S=$\verb|?...*o...o<| with possible zero \verb|o|'s between \verb|*| and \verb|<|.
\end{itemize}
As explained in (R) and (L), each time the right most element changes orientation, $r\circ M=M\circ r$. 

Since part a of the algorithm goes through all states exactly once with $k=l$ and $n=m$, we have that it also goes through all states exactly once with $k=l+1$ and $n=m+1$. The b-part of the algorithm follows from the a-part and the above theorem.
\end{proof}
\subsection{The tale of the wolves in the staircase}
The following description is equivalent to the rules in the algorithm.
Imagine $k$ wolves in a staircase with $n$ steps.
There is only one wolf on each step and each wolf can only face upwards or downwards.
The highest ranked wolf is the wolf that is located highest in the staircase at a present time.
This wolf is the only one that can initiate a move by itself.
If it cannot move, it must wait until another wolf has moved to an empty step.
A wolf can only go forwards to the nearest empty step.
However, it cannot pass a wolf that faces it.
It can pass as many wolves as it can as long as there is an empty step to go to and it passes all wolfs in between from behind.
Wolves that are passed by another wolf shall turn upwards.
A wolf that cannot move with these rules shall turn around and howl.
The bark signals that the closest wolf downstairs should try to move.
If the lowest ranked wolf barks, no wolves move.
The wolfs are now at the top of the staircase.
\subsection{Applying the algorithm on multisets}
Given the multiset $112233$ our algorithm gives the result as shown in \cref{tab:4}.
First start with the ones and treat the remaining elements as empty places.
For each time the algorithm exits, one step of the algorithm is applied on the string where the ones are removed and the twos are represented by \verb|>| or \verb|<|, the threes are treated as empty places.
\cref{tab:4} shows the process.
The table is read columnwise from the top to the bottom and from the left to the right.
\begin{table}[b!]
\begin{center}
\begin{tabular}{|c||c|c|c|c|c|c|}\hline
\verb|</>|&\multicolumn{6}{|c|}{States}\\\hline
2&\verb|>>33|&\verb|>3>3|&\verb|>33>|&\verb|3>3<|&\verb|3><3|&\verb|33>>|\\
\hline\hline
1&\verb|>>2233|&\verb|2323<<|&\verb|>>2332|&\verb|3232<<|&\verb|>>3223|&\verb|3322<<|\\
&\verb|>2>233|&\verb|232<>3|&\verb|>2>332|&\verb|323<>2|&\verb|>3>223|&\verb|332<>2|\\
&\verb|>22>33|&\verb|232<3>|&\verb|>23>32|&\verb|323<2>|&\verb|>32>23|&\verb|332<2>|\\
&\verb|>223>3|&\verb|23<23<|&\verb|>233>2|&\verb|32<32<|&\verb|>322>3|&\verb|33<22<|\\
&\verb|>2233>|&\verb|23<2<3|&\verb|>2332>|&\verb|32<3<2|&\verb|>3223>|&\verb|33<2<2|\\
&\verb|2>233<|&\verb|23<<23|&\verb|2>332<|&\verb|32<<32|&\verb|3>223<|&\verb|33<<22|\\
&\verb|2>23<3|&\verb|2<>323|&\verb|2>33<2|&\verb|3<>232|&\verb|3>22<3|&\verb|3<>322|\\
&\verb|2>2<33|&\verb|2<3>23|&\verb|2>3<32|&\verb|3<2>32|&\verb|3>2<23|&\verb|3<3>22|\\
&\verb|2><233|&\verb|2<32>3|&\verb|2><332|&\verb|3<23>2|&\verb|3><223|&\verb|3<32>2|\\
&\verb|22>>33|&\verb|2<323>|&\verb|23>>32|&\verb|3<232<|&\verb|32>>23|&\verb|3<322<|\\
&\verb|22>3>3|&\verb|<2323<|&\verb|23>3>2|&\verb|<3232<|&\verb|32>2>3|&\verb|<3322<|\\
&\verb|22>33>|&\verb|<232<3|&\verb|23>32>|&\verb|<323<2|&\verb|32>23>|&\verb|<332<2|\\
&\verb|223>3<|&\verb|<23<23|&\verb|233>2<|&\verb|<32<32|&\verb|322>3<|&\verb|<33<22|\\
&\verb|223><3|&\verb|<2<323|&\verb|233><2|&\verb|<3<232|&\verb|322><3|&\verb|<3<322|\\
&\verb|2233>>|&\verb|<<2323|&\verb|2332>>|&\verb|<<3232|&\verb|3223>>|&\verb|<<3322|\\\hline
\end{tabular}
\end{center}
\caption{The table shows how the algorithm permutes the multiset $1^22^23^2$.}
\label{tab:4}
\end{table}
The extension of our algorithm to arbitrary finite multisets is as following.
\begin{algorithm}
\label{alg:2}
Given the multiset $S=[1,\ldots,k]$ with $n$ elements and $S_{i}\leq S_{j}$ for $1\leq i<j\leq n$. The following algorithm traverses all permutations $P$ of $S$ exactly once.
\begin{enumerate}
\item Let $P=S$ and define the array $V=[1,\ldots,1]$ with $n$ elements $V_i=1$.
\item Set $T=1$. This is the active element type.
\label{A2_1}
\item Set $M=n+1$.
\label{A2_2}
\item Set $M:=\displaystyle\max_{i<M,P_i=T}i$.
\label{A2_3}
\item If $M$ does not exist:
\begin{itemize}
\item Set $T:=T+1$
\item If $T>k$. Exit the algorithm.
\item Go to step~(\ref{A2_2})
\end{itemize}
\item Let $N:=\begin{cases}
\displaystyle\min_{i>M,P_i>T}i&\mbox{if $V_M=1$}\\
\displaystyle\max_{i<M,P_i>T}i&\mbox{if $V_M=-1$}
\end{cases}$
\item If $N$ does not exist, go to step~(\ref{A2_4}).
\item If $V_M=V_i$ for all $i$ between $M$ and $N$
\begin{itemize}
\item Swap $P_M$ with $P_N$ and $V_M$ with $V_N$
\item Set $V_i=1$ for all $i$ between $M$ and $N$.
\item Go to step~(\ref{A2_1})
\end{itemize}
\item Set $V_M:=-V_M$.
\label{A2_4}
\item Go to step~(\ref{A2_3})
\end{enumerate}
\end{algorithm}
\subsection{A circular Gray code for a multiset.}
Our algorithm is circular in the case $m_{k-1}=m_k=1$.
\cref{fig:GrayCode} shows the result of the algorithm applied on the multiset $1^22^234$.
\section{Generating the polynomials}
\label{sec:gen_pol}
In this section, we apply the algorithm on the Young tableaux $B_{2,2,2,2}$.
Apply the algorithm on each column.
In the first column, the algorithm identifies 1 with 2, and 5 with 6.
The algorithm gives the sequence: 
$$1256,1526,1562,5162,5126,5621.$$
Similarly, for the second column of $B_{2,2,2,2}$, it produces the list:
$$3478,3748,3784,7384,7348,7843.$$
\cref{tab:5} shows the result.
\newcommand{\young}[8]{\begin{tikzpicture}[
baseline=($(A.base)!.5!(H.base)$),
auto,node distance=4mm,semithick,inner sep=2pt,bend angle=45]
\node [box] (A)              {$#1$};
\node [box] (B) [below of=A] {$#2$};
\node [box] (E) [below of=B] {$#3$};
\node [box] (F) [below of=E] {$#4$};
\node [box] (C) [right of=A] {$#5$};
\node [box] (D) [below of=C] {$#6$};
\node [box] (G) [below of=D] {$#7$};
\node [box] (H) [below of=G] {$#8$};
\end{tikzpicture}
}
\begin{table}[b!]
\begin{align*}
\phantom{+}\young12563478
-\young15263478
+\young15623478
-\young51623478
+\young51263478
-\young56213478\\
+\young56213748
-\young51263748
+\young51623748
-\young15623748
+\young15263748
-\young12563748\\
+\young12563784
-\young15263784
+\young15623784
-\young51623784
+\young51263784
-\young56213784\\
+\young56217384
-\young51267384
+\young51627384
-\young15627384
+\young15267384
-\young12567384\\
+\young12567348
-\young15267348
+\young15627348
-\young51627348
+\young51267348
-\young56217348\\
+\young56217843
-\young51267843
+\young51627843
-\young15627843
+\young15267843
-\young12567843\\
\end{align*}
\caption{The multi-set permutations of the diagram 
$B_{2,2,2,2}$.
The multi-sets are defined by grouping the elements in the equivalence classes
$\{1,2\}$,
$\{3,4\}$,
$\{5,6\}$, and
$\{7,8\}$.
The signs are the signatures of the permutations.
}
\label{tab:5}
\end{table}
This gives the polynomial
\begin{dmath*}
P(R)
=R_{1212}R_{3434}
-R_{1312}R_{2434}
+R_{1412}R_{2334}
-R_{2412}R_{1334}
+R_{2312}R_{1434}
-R_{4312}R_{1234}
+R_{4313}R_{1224}
-R_{2313}R_{1424}
+R_{2413}R_{1324}
-R_{1413}R_{2324}
+R_{1313}R_{2424}
-R_{1213}R_{3424}
+R_{1214}R_{3423}
-R_{1314}R_{2423}
+R_{1414}R_{2323}
-R_{2414}R_{1323}
+R_{2314}R_{1423}
-R_{4314}R_{1223}
+R_{4324}R_{1213}
-R_{2324}R_{1413}
+R_{2424}R_{1313}
-R_{1424}R_{2313}
+R_{1324}R_{2413}
-R_{1224}R_{3413}
+R_{1223}R_{3414}
-R_{1323}R_{2414}
+R_{1423}R_{2314}
-R_{2423}R_{1314}
+R_{2323}R_{1414}
-R_{4323}R_{1214}
+R_{4343}R_{1212}.
-R_{2343}R_{1412}
+R_{2443}R_{1312}
-R_{1443}R_{2312}
+R_{1343}R_{2412}
-R_{1243}R_{3412}
\end{dmath*}
By using the symmetries of $R$, this simplifies to
\begin{dmath*}
P(R)/2
={R_{1234}}^2
+{R_{1423}}^2
+{R_{1324}}^2
+R_{1212}R_{3434}
+R_{1313}R_{2424}
+R_{1414}R_{2323}
+2R_{1214}R_{2334}
+2R_{1223}R_{1434}
-2R_{1213}R_{2434}
-2R_{1224}R_{1334}
-2R_{1314}R_{2324}
-2R_{1323}R_{1424}.
\end{dmath*}
This is equal to the polynomial $I_{B_{2,2,2,2}}$ on page 118 in~\cite{agaoka85}.
For the larger Young tableau,
$$
B_{5,5,5,5,4,4}=
\begin{tikzpicture}[
baseline=($(a1.base)!.5!(f1.base)$),
auto,node distance=4.8mm,semithick,inner sep=2pt,bend angle=45]
\node [box] (a1) {1};
\node [box] (a2) [right of=a1] {3};
\node [box] (a3) [right of=a2] {5};
\node [box] (a4) [right of=a3] {7};
\node [box] (a5) [right of=a4] {9};
\node [box] (b1) [below of=a1] {2};
\node [box] (b2) [right of=b1] {4};
\node [box] (b3) [right of=b2] {6};
\node [box] (b4) [right of=b3] {8};
\node [box] (b5) [right of=b4] {10};
\node [box] (c1) [below of=b1] {11};
\node [box] (c2) [right of=c1] {13};
\node [box] (c3) [right of=c2] {15};
\node [box] (c4) [right of=c3] {17};
\node [box] (c5) [right of=c4] {19};
\node [box] (d1) [below of=c1] {12};
\node [box] (d2) [right of=d1] {14};
\node [box] (d3) [right of=d2] {16};
\node [box] (d4) [right of=d3] {18};
\node [box] (d5) [right of=d4] {20};
\node [box] (e1) [below of=d1] {21};
\node [box] (e2) [right of=e1] {23};
\node [box] (e3) [right of=e2] {25};
\node [box] (e4) [right of=e3] {27};
\node [box] (f1) [below of=e1] {22};
\node [box] (f2) [right of=f1] {24};
\node [box] (f3) [right of=f2] {26};
\node [box] (f4) [right of=f3] {28};
\end{tikzpicture}
$$
the saved amount of calculation is significant. The size of $\mathfrak{V}_{B_{5,5,5,5,4,4}}$ is $(6!)^24!=6{,}449{,}725{,}440{,}000$.
By using the first two symmetries in the curvature tensor, the algorithm reduces the complexity to $\frac{(6!)^24!}{2^{14}}=393{,}660{,}000$.
Using the remaining symmetries will reduce the complexity of the calculations even more.
That is a topic for another paper.
\section{Discussion}
\subsection{Comparison to Eades McKay}
Our algorithm is quite similar to Eades-McKay's algorithm.
Both methods uses transpositions $(i,j)$ only.
Eades-McKay's algorithm require the homogeneous transpositions condition in~\cite{Ruskey2003}.
Our algorithm follows a strong version of the homogeneous transpositions condition for the empty spaces.
The elements between two interchanging elements must be equal to the smallest of the two interchanging elements.
\subsection{Notes on an algorithm of Takaoka}
Takaoka's algorithm~\cite{Takaoka1999} for generating all combinations of $k$ elements in a finite set seems to be a special case (a.) of our \cref{alg:1}.
Although, Takaoka's algorithm is for generating combinations, both his and our algorithms produce the same result for the multisets $1^20^4$ and $1^30^3$.
It remains to prove that those algorithms are the same.
Takaoka's algorithm~\cite{Takaoka1999_a} for producing all permutations of a multiset is different from our algorithm.
\subsection{Inspiration from Steinhaus-Johnson-Trotter}
\begin{figure}[htb!]
\begin{center}
\begin{tikzpicture}[scale=0.9,auto,node distance=1.2cm,inner sep=1pt,bend angle=45]
\node [my state] (Ao) {1234};
\node [my state] (Bo) at ([shift=({  60:5cm})]Ao) {2134};
\node [my state] (Co) at ([shift=({   0:5cm})]Bo) {2314};
\node [my state] (Do) at ([shift=({ -60:5cm})]Co) {3214};
\node [my state] (Eo) at ([shift=({-120:5cm})]Do) {3124};
\node [my state] (Fo) at ([shift=({ 180:5cm})]Eo) {1324};
\node [my state] (Ai) at ([shift=({ 30:1.5cm})]Ao) {1243};
\node [my state] (Bi) at ([shift=({270:1.5cm})]Bo) {2143};
\node [my state] (Ci) at ([shift=({270:1.5cm})]Co) {2341};
\node [my state] (Di) at ([shift=({150:1.5cm})]Do) {3241};
\node [my state] (Ei) at ([shift=({150:1.5cm})]Eo) {3142};
\node [my state] (Fi) at ([shift=({ 30:1.5cm})]Fo) {1342};
\node [my state] (Ak) at ([shift=({300:1.5cm})]Ai) {1423};
\node [my state] (Bk) at ([shift=({  0:1.5cm})]Bi) {2413};
\node [my state] (Ck) at ([shift=({180:1.5cm})]Ci) {2431};
\node [my state] (Dk) at ([shift=({240:1.5cm})]Di) {3421};
\node [my state] (Ek) at ([shift=({ 60:1.5cm})]Ei) {3412};
\node [my state] (Fk) at ([shift=({120:1.5cm})]Fi) {1432};
\node [my state] (Al) at ([shift=({ 30:1.3cm})]Ak) {4123};
\node [my state] (Bl) at ([shift=({270:1.3cm})]Bk) {4213};
\node [my state] (Cl) at ([shift=({270:1.3cm})]Ck) {4231};
\node [my state] (Dl) at ([shift=({150:1.3cm})]Dk) {4321};
\node [my state] (El) at ([shift=({150:1.3cm})]Ek) {4312};
\node [my state] (Fl) at ([shift=({ 30:1.3cm})]Fk) {4132};
\path (Ao) edge node {(1,2)} (Bo)
      (Bo) edge node {(2,3)} (Co)
      (Co) edge node {(1,2)} (Do)
      (Do) edge node {(2,3)} (Eo)
      (Eo) edge node {(1,2)} (Fo)
      (Fo) edge node {(2,3)} (Ao);
\path (Ai) edge node {(3,4)} (Ao)
      (Bo) edge node {(3,4)} (Bi)
      (Ci) edge node {(3,4)} (Co)
      (Do) edge node {(3,4)} (Di)
      (Ei) edge node {(3,4)} (Eo)
      (Fo) edge node {(3,4)} (Fi);
\path (Bi) edge node {(1,2)} (Ai)
      (Di) edge node {(1,2)} (Ci)
      (Fi) edge node {(1,2)} (Ei);
\path (Ai) edge node {(2,3)} (Ak)
      (Bk) edge node {(2,3)} (Bi)
      (Ci) edge node {(2,3)} (Ck)
      (Dk) edge node {(2,3)} (Di)
      (Ei) edge node {(2,3)} (Ek)
      (Fk) edge node {(2,3)} (Fi);
\path (Bk) edge node {(3,4)} (Ck)
      (Dk) edge node {(3,4)} (Ek)
      (Fk) edge node {(3,4)} (Ak);
\path (Al) edge node {(1,2)} (Ak)
      (Bk) edge node {(1,2)} (Bl)
      (Cl) edge node {(1,2)} (Ck)
      (Dk) edge node {(1,2)} (Dl)
      (El) edge node {(1,2)} (Ek)
      (Fk) edge node {(1,2)} (Fl);
\path (Bl) edge node {(2,3)} (Al)
      (Cl) edge node {(3,4)} (Bl)
      (Dl) edge node {(2,3)} (Cl)
      (El) edge node {(3,4)} (Dl)
      (Fl) edge node {(2,3)} (El)
      (Al) edge node {(3,4)} (Fl);
\draw[blue,ultra thick,->] (Ao) -- (Bo) -- (Co) -- (Ci) -- (Di) -- (Do) -- (Eo) -- (Fo) -- (Fi) -- (Ei) -- (Ek) -- (Dk) -- (Dl) -- (El) -- (Fl) -- (Fk) -- (Ak) -- (Al) -- (Bl) -- (Cl) -- (Ck) -- (Bk) -- (Bi) -- (Ai) -- (Ao);
\end{tikzpicture}
\end{center}
\caption{The vertices of the graph are all permutations of 1234. The edges are transpositions. Our algorithm is marked with a thick line which forms a Hamilton circuit. This circuit is exactly the same as the one produced by the Steinhaus-Johnson-Trotter algorithm~\cite{Selmer_M_Johnson:1963,Trotter:1962}.
The graph is homomorphic to the truncated octahedron, one of Archimedes solids as shown in \cref{fig:TrunkOcth}.
}
\label{fig:1234}
\end{figure}
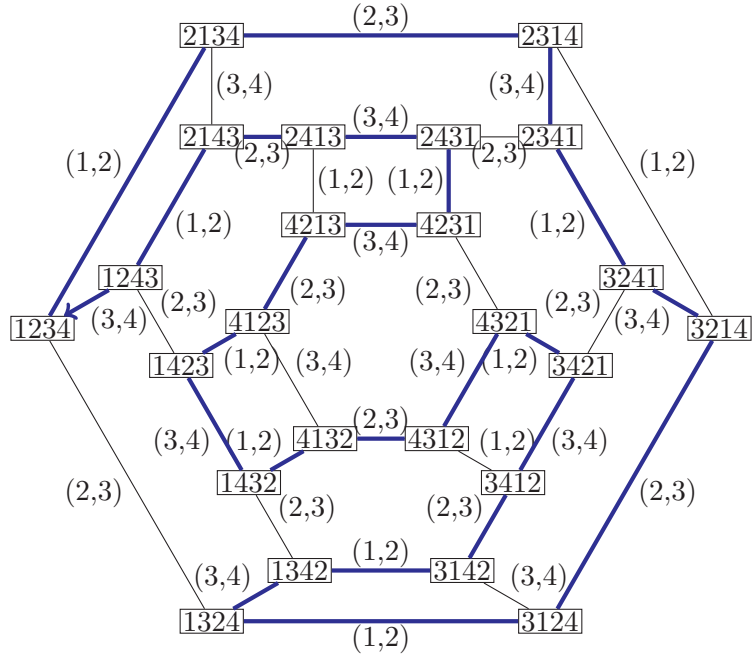

Our method was inspired by the Steinhaus-Johnson-Trotter algorithm~\cite{Selmer_M_Johnson:1963,Trotter:1962}.
In their algorithm, the element $1$ goes back and forth. Each time it turns around, the next element $2$ takes one step and so on.
\cref{fig:1234} shows that their algorithm and our algorithm are equal when applied on the set $\{1,2,3,4\}$.
We believe that our algorithm is a generalisation of their algorithm.
\subsection{Unanswered question}
Each step in our algorithm is a transposition $(i,j)$.
In other words, the Hamming distance is two
between successive elements in the list $B$ of all permutations of the multiset $1^k0^{n-k}$.
This is the first order closeness measure.
The width of a transposition $(i,j)$ is defined as $w_{(i,j)}=|i-j|$.
The sum of the widths for all transitions in the algorithm is called the total motion
$W_B=\sum_{\sigma}w_{\sigma}$.
It is an open question if the presented algorithm has the minimal total motion.
In the case $n=6$ and $k=3$, the Eades–McKay algorithm has 16 adjacent transpositions, two moves of distance two and one move of distance 3.
The total distance moved is therefore $16+2\times2+3=23$.
Our method has 15 adjacent transpositions and 4 moves of distance 2.
The total distance moved is for our method is $15+4\times2=23$.
Another condition is to only accept adjacent or transpositions with width 2.
A known algorithm~\cite{Ruskey2003} that have this condition has total distance moved 25, for $n=6$ and $k=3$.

Let $B(n,k)$ be the list of bit-strings produced by our method.
Experiments show that:
\begin{itemize}
\item[a)]
$W_{E(n,n-k)}=W_{B(n,k)}$ for all $n<30$ and $0<k<n$.
\item[b)]
$W_{E(n,k)}>W_{B(n,k)}$ for all $n<30$ and $0<k<n/2$.
\item[c)]
$W_{E(n,k)}<W_{B(n,k)}$ for all $n<30$ and $n/2<k<n$.
\end{itemize}
It is not known which method has the lowest value for $W_{E(n,k)}$.
The values in the list above is limited to the case where $n<30$ and for two different algorithms only.

\begin{figure}[htb!]
\begin{center}
\tdplotsetmaincoords{80}{115}
\begin{tikzpicture}[scale=1.5,tdplot_main_coords]
\node (BC) at (-2,1,0) {1234};
\node (BE) at (-2,0,1) {2134};
\node (EB) at (-1,0,2) {2314};
\node (ED) at (0,-1,2) {2341};
\node (EA) at (1,0,2)  {3241};
\node (EC) at (0,1,2)  {3214};
\node (CE) at (0,2,1)  {3124};
\node (CB) at (-1,2,0) {1324};
\node (CF) at (0,2,-1) {1342};
\node (CA) at (1,2,0)  {3142};
\node (AC) at (2, 1,0) {3412};
\node (AE) at (2,0,1)  {3421};
\node (AD) at (2,-1,0) {4321};
\node (AF) at (2,0,-1) {4312};
\node (FA) at (1,0,-2) {4132};
\node (FC) at (0,1,-2) {1432};
\node (FB) at (-1,0,-2){1423};
\node (FD) at (0,-1,-2){4123};
\node (DF) at (0,-2,-1){4213};
\node (DA) at (1,-2,0) {4231};
\node (DE) at (0,-2,1) {2431};
\node (DB) at (-1,-2,0){2413};
\node (BD) at (-2,-1,0){2143};
\node (BF) at (-2,0,-1){1243};

\draw[color=blue,dashed,ultra thick](EA) -- (EC);
\draw[ultra thick]            (EC) -- node[above right] {$(2,3)$} (CE) ;
\draw[color=red,ultra thick]        (CE) -- node[right] {$(1,2)$} (CB);
\draw[color=blue,ultra thick]       (CB) -- node[right] {$(3,4)$} (CF);

\draw[color=red,dashed,ultra thick] (CF) -- (CA);
\draw[dashed,ultra thick]           (CA) -- (AC);
\draw[color=blue,dashed,ultra thick](AC) -- (AE);
\draw[color=red,dashed,ultra thick] (AE) -- (AD);
\draw[color=blue,dashed,ultra thick](AD) -- (AF);
\draw[dashed,ultra thick]           (AF) -- (FA);
\draw[color=red,ultra thick]        (FA) -- (FC);
\draw[color=blue,ultra thick]       (FC) -- (FB);
\draw[color=red,ultra thick]        (FB) -- (FD);
\draw[ultra thick]                  (FD) -- (DF);
\draw[color=blue,ultra thick]       (DF) -- (DA);
\draw[color=red,ultra thick]        (DA) -- (DE);
\draw[color=blue,ultra thick]       (DE) -- (DB);
\draw[ultra thick]                  (DB) -- (BD);
\draw[color=red,ultra thick]        (BD) -- (BF);
\draw[color=blue,ultra thick,->]    (BF) -- (BC);
\draw[color=red,ultra thick]        (BC) -- (BE);
\draw[ultra thick]                  (BE) -- (EB);
\draw[color=blue,ultra thick]       (EB) -- (ED);
\draw[color=red,dashed,ultra thick] (ED) -- (EA);

\draw[color=red] (DF) -- (DB);
\draw[color=red] (EC) -- (EB);
\draw[color=blue] (FD) -- (FA);
\draw[color=blue,dashed] (CA) -- (CE);
\draw[color=red,dashed](AF) -- (AC);
\draw[color=blue] (BE) -- (BD);
\draw[dashed] (AD) -- (DA);
\draw[dashed] (AE) -- (EA);
\draw (BC) -- (CB);
\draw (BF) -- (FB);
\draw (CF) -- (FC);
\draw (DE) -- (ED);
\end{tikzpicture}
\end{center}
\caption{The graph in \cref{fig:1234} is drawn as a truncated octahedron. The edges are adjacent transpositions. Our algorithm is marked with a thick line which forms a Hamilton circuit.
}
\label{fig:TrunkOcth}
\end{figure}
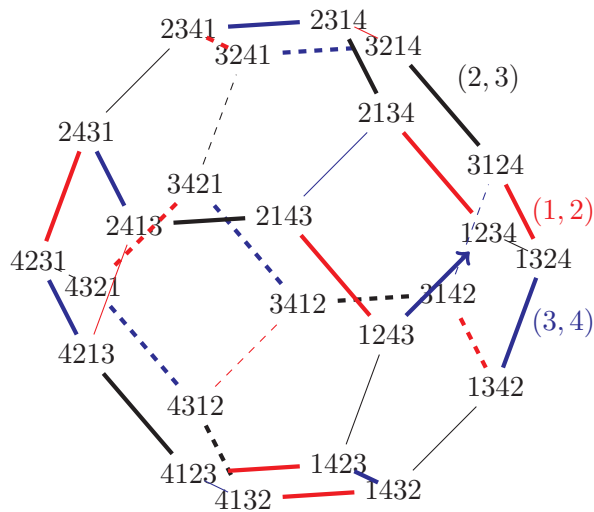

\def\sector#1#2#3#4#5{
\fill[#3] (xy polar cs:angle=#1,radius=#4) --  (xy polar cs:angle=#1,radius=#5) arc (#1:#2:#5) --  (xy polar cs:angle=#2,radius=#4) arc (#2:#1:#4);
}
\def\msector#1#2#3#4#5#6#7#8{
\sector{#1}{#2}{#3}{2.0}{2.5}
\sector{#1}{#2}{#4}{2.5}{3.0}
\sector{#1}{#2}{#5}{3.0}{3.5}
\sector{#1}{#2}{#6}{3.5}{4.0}
\sector{#1}{#2}{#7}{4.0}{4.5}
\sector{#1}{#2}{#8}{4.5}{5.0}
\draw[thin,dashed] (xy polar cs:angle=#1,radius=2) -- (xy polar cs:angle=#1,radius=5) arc (#1:#2:5) -- (xy polar cs:angle=#2,radius=2) arc (#2:#1:2);
}

\def\msectorii#1#2#3#4#5#6{
\sector{#1}{#2}{#3}{0.5}{0.75}
\sector{#1}{#2}{#4}{0.75}{1.0}
\sector{#1}{#2}{#5}{1.0}{1.25}
\sector{#1}{#2}{#6}{1.25}{1.5}
\draw[thin,dashed] (xy polar cs:angle=#1,radius=.5) -- (xy polar cs:angle=#1,radius=1.5) arc (#1:#2:1.5) -- (xy polar cs:angle=#2,radius=.5) arc (#2:#1:.5);
}

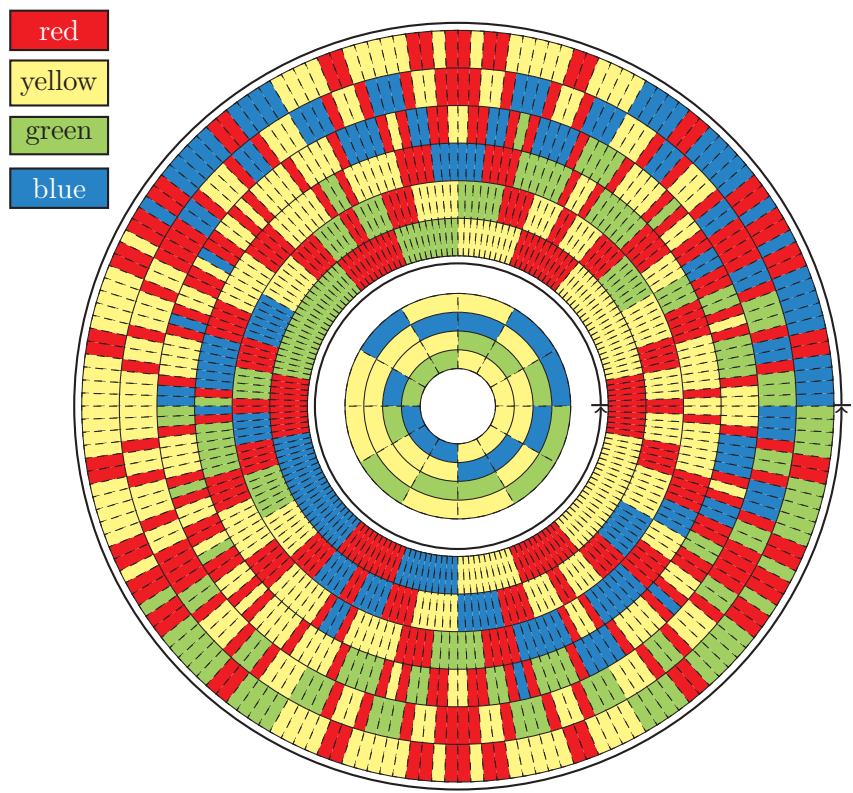
\begin{figure}[h!]
\begin{center}
\begin{tikzpicture}
\path [draw,|->,thick] (5.1,0) arc (0:360:5.1);
\path [draw,|->,thick] (1.9,0) arc (0:360:1.9);
\msector{  0}{  2}{R}{R}{Y}{Y}{G}{B}
\msector{  2}{  4}{R}{Y}{R}{Y}{G}{B}
\msector{  4}{  6}{R}{Y}{Y}{R}{G}{B}
\msector{  6}{  8}{R}{Y}{Y}{G}{R}{B}
\msector{  8}{ 10}{R}{Y}{Y}{G}{B}{R}
\msector{ 10}{ 12}{Y}{R}{Y}{G}{B}{R}
\msector{ 12}{ 14}{Y}{R}{Y}{G}{R}{B}
\msector{ 14}{ 16}{Y}{R}{Y}{R}{G}{B}
\msector{ 16}{ 18}{Y}{R}{R}{Y}{G}{B}
\msector{ 18}{ 20}{Y}{Y}{R}{R}{G}{B}
\msector{ 20}{ 22}{Y}{Y}{R}{G}{R}{B}
\msector{ 22}{ 24}{Y}{Y}{R}{G}{B}{R}
\msector{ 24}{ 26}{Y}{Y}{G}{R}{B}{R}
\msector{ 26}{ 28}{Y}{Y}{G}{R}{R}{B}
\msector{ 28}{ 30}{Y}{Y}{G}{B}{R}{R}
\msector{ 30}{ 32}{Y}{G}{Y}{B}{R}{R}
\msector{ 32}{ 34}{Y}{G}{Y}{R}{R}{B}
\msector{ 34}{ 36}{Y}{G}{Y}{R}{B}{R}
\msector{ 36}{ 38}{Y}{G}{R}{Y}{B}{R}
\msector{ 38}{ 40}{Y}{G}{R}{Y}{R}{B}
\msector{ 40}{ 42}{Y}{G}{R}{R}{Y}{B}
\msector{ 42}{ 44}{Y}{R}{R}{G}{Y}{B}
\msector{ 44}{ 46}{Y}{R}{G}{R}{Y}{B}
\msector{ 46}{ 48}{Y}{R}{G}{Y}{R}{B}
\msector{ 48}{ 50}{Y}{R}{G}{Y}{B}{R}
\msector{ 50}{ 52}{R}{Y}{G}{Y}{B}{R}
\msector{ 52}{ 54}{R}{Y}{G}{Y}{R}{B}
\msector{ 54}{ 56}{R}{Y}{G}{R}{Y}{B}
\msector{ 56}{ 58}{R}{Y}{R}{G}{Y}{B}
\msector{ 58}{ 60}{R}{R}{Y}{G}{Y}{B}
\msector{ 60}{ 62}{R}{R}{Y}{G}{B}{Y}
\msector{ 62}{ 64}{R}{Y}{R}{G}{B}{Y}
\msector{ 64}{ 66}{R}{Y}{G}{R}{B}{Y}
\msector{ 66}{ 68}{R}{Y}{G}{B}{R}{Y}
\msector{ 68}{ 70}{R}{Y}{G}{B}{Y}{R}
\msector{ 70}{ 72}{Y}{R}{G}{B}{Y}{R}
\msector{ 72}{ 74}{Y}{R}{G}{B}{R}{Y}
\msector{ 74}{ 76}{Y}{R}{G}{R}{B}{Y}
\msector{ 76}{ 78}{Y}{R}{R}{G}{B}{Y}
\msector{ 78}{ 80}{Y}{G}{R}{R}{B}{Y}
\msector{ 80}{ 82}{Y}{G}{R}{B}{R}{Y}
\msector{ 82}{ 84}{Y}{G}{R}{B}{Y}{R}
\msector{ 84}{ 86}{Y}{G}{B}{R}{Y}{R}
\msector{ 86}{ 88}{Y}{G}{B}{R}{R}{Y}
\msector{ 88}{ 90}{Y}{G}{B}{Y}{R}{R}
\msector{ 90}{ 92}{G}{Y}{B}{Y}{R}{R}
\msector{ 92}{ 94}{G}{Y}{B}{R}{R}{Y}
\msector{ 94}{ 96}{G}{Y}{B}{R}{Y}{R}
\msector{ 96}{ 98}{G}{Y}{R}{B}{Y}{R}
\msector{ 98}{100}{G}{Y}{R}{B}{R}{Y}
\msector{100}{102}{G}{Y}{R}{R}{B}{Y}
\msector{102}{104}{G}{R}{R}{Y}{B}{Y}
\msector{104}{106}{G}{R}{Y}{R}{B}{Y}
\msector{106}{108}{G}{R}{Y}{B}{R}{Y}
\msector{108}{110}{G}{R}{Y}{B}{Y}{R}
\msector{110}{112}{R}{G}{Y}{B}{Y}{R}
\msector{112}{114}{R}{G}{Y}{B}{R}{Y}
\msector{114}{116}{R}{G}{Y}{R}{B}{Y}
\msector{116}{118}{R}{G}{R}{Y}{B}{Y}
\msector{118}{120}{R}{R}{G}{Y}{B}{Y}
\msector{120}{122}{R}{R}{G}{Y}{Y}{B}
\msector{122}{124}{R}{G}{R}{Y}{Y}{B}
\msector{124}{126}{R}{G}{Y}{R}{Y}{B}
\msector{126}{128}{R}{G}{Y}{Y}{R}{B}
\msector{128}{130}{R}{G}{Y}{Y}{B}{R}
\msector{130}{132}{G}{R}{Y}{Y}{B}{R}
\msector{132}{134}{G}{R}{Y}{Y}{R}{B}
\msector{134}{136}{G}{R}{Y}{R}{Y}{B}
\msector{136}{138}{G}{R}{R}{Y}{Y}{B}
\msector{138}{140}{G}{Y}{R}{R}{Y}{B}
\msector{140}{142}{G}{Y}{R}{Y}{R}{B}
\msector{142}{144}{G}{Y}{R}{Y}{B}{R}
\msector{144}{146}{G}{Y}{Y}{R}{B}{R}
\msector{146}{148}{G}{Y}{Y}{R}{R}{B}
\msector{148}{150}{G}{Y}{Y}{B}{R}{R}
\msector{150}{152}{G}{B}{Y}{Y}{R}{R}
\msector{152}{154}{G}{B}{Y}{R}{R}{Y}
\msector{154}{156}{G}{B}{Y}{R}{Y}{R}
\msector{156}{158}{G}{B}{R}{Y}{Y}{R}
\msector{158}{160}{G}{B}{R}{Y}{R}{Y}
\msector{160}{162}{G}{B}{R}{R}{Y}{Y}
\msector{162}{164}{G}{R}{R}{B}{Y}{Y}
\msector{164}{166}{G}{R}{B}{R}{Y}{Y}
\msector{166}{168}{G}{R}{B}{Y}{R}{Y}
\msector{168}{170}{G}{R}{B}{Y}{Y}{R}
\msector{170}{172}{R}{G}{B}{Y}{Y}{R}
\msector{172}{174}{R}{G}{B}{Y}{R}{Y}
\msector{174}{176}{R}{G}{B}{R}{Y}{Y}
\msector{176}{178}{R}{G}{R}{B}{Y}{Y}
\msector{178}{180}{R}{R}{G}{B}{Y}{Y}
\msector{180}{182}{R}{R}{B}{G}{Y}{Y}
\msector{182}{184}{R}{B}{R}{G}{Y}{Y}
\msector{184}{186}{R}{B}{G}{R}{Y}{Y}
\msector{186}{188}{R}{B}{G}{Y}{R}{Y}
\msector{188}{190}{R}{B}{G}{Y}{Y}{R}
\msector{190}{192}{B}{R}{G}{Y}{Y}{R}
\msector{192}{194}{B}{R}{G}{Y}{R}{Y}
\msector{194}{196}{B}{R}{G}{R}{Y}{Y}
\msector{196}{198}{B}{R}{R}{G}{Y}{Y}
\msector{198}{200}{B}{G}{R}{R}{Y}{Y}
\msector{200}{202}{B}{G}{R}{Y}{R}{Y}
\msector{202}{204}{B}{G}{R}{Y}{Y}{R}
\msector{204}{206}{B}{G}{Y}{R}{Y}{R}
\msector{206}{208}{B}{G}{Y}{R}{R}{Y}
\msector{208}{210}{B}{G}{Y}{Y}{R}{R}
\msector{210}{212}{B}{Y}{Y}{G}{R}{R}
\msector{212}{214}{B}{Y}{Y}{R}{R}{G}
\msector{214}{216}{B}{Y}{Y}{R}{G}{R}
\msector{216}{218}{B}{Y}{R}{Y}{G}{R}
\msector{218}{220}{B}{Y}{R}{Y}{R}{G}
\msector{220}{222}{B}{Y}{R}{R}{Y}{G}
\msector{222}{224}{B}{R}{R}{Y}{Y}{G}
\msector{224}{226}{B}{R}{Y}{R}{Y}{G}
\msector{226}{228}{B}{R}{Y}{Y}{R}{G}
\msector{228}{230}{B}{R}{Y}{Y}{G}{R}
\msector{230}{232}{R}{B}{Y}{Y}{G}{R}
\msector{232}{234}{R}{B}{Y}{Y}{R}{G}
\msector{234}{236}{R}{B}{Y}{R}{Y}{G}
\msector{236}{238}{R}{B}{R}{Y}{Y}{G}
\msector{238}{240}{R}{R}{B}{Y}{Y}{G}
\msector{240}{242}{R}{R}{B}{Y}{G}{Y}
\msector{242}{244}{R}{B}{R}{Y}{G}{Y}
\msector{244}{246}{R}{B}{Y}{R}{G}{Y}
\msector{246}{248}{R}{B}{Y}{G}{R}{Y}
\msector{248}{250}{R}{B}{Y}{G}{Y}{R}
\msector{250}{252}{B}{R}{Y}{G}{Y}{R}
\msector{252}{254}{B}{R}{Y}{G}{R}{Y}
\msector{254}{256}{B}{R}{Y}{R}{G}{Y}
\msector{256}{258}{B}{R}{R}{Y}{G}{Y}
\msector{258}{260}{B}{Y}{R}{R}{G}{Y}
\msector{260}{262}{B}{Y}{R}{G}{R}{Y}
\msector{262}{264}{B}{Y}{R}{G}{Y}{R}
\msector{264}{266}{B}{Y}{G}{R}{Y}{R}
\msector{266}{268}{B}{Y}{G}{R}{R}{Y}
\msector{268}{270}{B}{Y}{G}{Y}{R}{R}
\msector{270}{272}{Y}{B}{G}{Y}{R}{R}
\msector{272}{274}{Y}{B}{G}{R}{R}{Y}
\msector{274}{276}{Y}{B}{G}{R}{Y}{R}
\msector{276}{278}{Y}{B}{R}{G}{Y}{R}
\msector{278}{280}{Y}{B}{R}{G}{R}{Y}
\msector{280}{282}{Y}{B}{R}{R}{G}{Y}
\msector{282}{284}{Y}{R}{R}{B}{G}{Y}
\msector{284}{286}{Y}{R}{B}{R}{G}{Y}
\msector{286}{288}{Y}{R}{B}{G}{R}{Y}
\msector{288}{290}{Y}{R}{B}{G}{Y}{R}
\msector{290}{292}{R}{Y}{B}{G}{Y}{R}
\msector{292}{294}{R}{Y}{B}{G}{R}{Y}
\msector{294}{296}{R}{Y}{B}{R}{G}{Y}
\msector{296}{298}{R}{Y}{R}{B}{G}{Y}
\msector{298}{300}{R}{R}{Y}{B}{G}{Y}
\msector{300}{302}{R}{R}{Y}{B}{Y}{G}
\msector{302}{304}{R}{Y}{R}{B}{Y}{G}
\msector{304}{306}{R}{Y}{B}{R}{Y}{G}
\msector{306}{308}{R}{Y}{B}{Y}{R}{G}
\msector{308}{310}{R}{Y}{B}{Y}{G}{R}
\msector{310}{312}{Y}{R}{B}{Y}{G}{R}
\msector{312}{314}{Y}{R}{B}{Y}{R}{G}
\msector{314}{316}{Y}{R}{B}{R}{Y}{G}
\msector{316}{318}{Y}{R}{R}{B}{Y}{G}
\msector{318}{320}{Y}{B}{R}{R}{Y}{G}
\msector{320}{322}{Y}{B}{R}{Y}{R}{G}
\msector{322}{324}{Y}{B}{R}{Y}{G}{R}
\msector{324}{326}{Y}{B}{Y}{R}{G}{R}
\msector{326}{328}{Y}{B}{Y}{R}{R}{G}
\msector{328}{330}{Y}{B}{Y}{G}{R}{R}
\msector{330}{332}{Y}{Y}{B}{G}{R}{R}
\msector{332}{334}{Y}{Y}{B}{R}{R}{G}
\msector{334}{336}{Y}{Y}{B}{R}{G}{R}
\msector{336}{338}{Y}{Y}{R}{B}{G}{R}
\msector{338}{340}{Y}{Y}{R}{B}{R}{G}
\msector{340}{342}{Y}{Y}{R}{R}{B}{G}
\msector{342}{344}{Y}{R}{R}{Y}{B}{G}
\msector{344}{346}{Y}{R}{Y}{R}{B}{G}
\msector{346}{348}{Y}{R}{Y}{B}{R}{G}
\msector{348}{350}{Y}{R}{Y}{B}{G}{R}
\msector{350}{352}{R}{Y}{Y}{B}{G}{R}
\msector{352}{354}{R}{Y}{Y}{B}{R}{G}
\msector{354}{356}{R}{Y}{Y}{R}{B}{G}
\msector{356}{358}{R}{Y}{R}{Y}{B}{G}
\msector{358}{360}{R}{R}{Y}{Y}{B}{G} \msectorii{0}{30}{Y}{Y}{G}{B}
\msectorii{30}{60}{Y}{G}{Y}{B}
\msectorii{60}{90}{Y}{G}{B}{Y}
\msectorii{90}{120}{G}{Y}{B}{Y}
\msectorii{120}{150}{G}{Y}{Y}{B}
\msectorii{150}{180}{G}{B}{Y}{Y}
\msectorii{180}{210}{B}{G}{Y}{Y}
\msectorii{210}{240}{B}{Y}{Y}{G}
\msectorii{240}{270}{B}{Y}{G}{Y}
\msectorii{270}{300}{Y}{B}{G}{Y}
\msectorii{300}{330}{Y}{B}{Y}{G}
\msectorii{330}{360}{Y}{Y}{B}{G}
\draw[thin] (0,0) circle (0.50cm);
\draw[thin] (0,0) circle (0.75cm);
\draw[thin] (0,0) circle (1.00cm);
\draw[thin] (0,0) circle (1.25cm);
\draw[thin] (0,0) circle (1.50cm);
\draw[thin] (0,0) circle (2.0cm);
\draw[thin] (0,0) circle (2.5cm);
\draw[thin] (0,0) circle (3.0cm);
\draw[thin] (0,0) circle (3.5cm);
\draw[thin] (0,0) circle (4.0cm);
\draw[thin] (0,0) circle (4.5cm);
\draw[thin] (0,0) circle (5.0cm);
\node (rect) at (-5.3,5.0) [draw,thick,minimum width=13mm,minimum height=4mm,fill=R] {\color{white}red};
\node (rect) at (-5.3,4.3) [draw,thick,minimum width=13mm,minimum height=4mm,fill=Y] {yellow};
\node (rect) at (-5.3,3.6) [draw,thick,minimum width=13mm,minimum height=4mm,fill=G] {green};
\node (rect) at (-5.3,2.9) [draw,thick,minimum width=13mm,minimum height=4mm,fill=B] {\color{white}blue};

\end{tikzpicture}
\end{center}
\caption{The picture shows the result of our algorithm applied on the multiset $112234$. Notice that the pattern is a circular Gray code. The colour codes are
1=red, 2=yellow, 3=green, and 4=blue. The inner circle shows the output from the algorithm for the multiset $1123$. Notice that the order of yellow, green, and blue in the inner figure is the same as the order of the same colors in the outer figure.}
\label{fig:GrayCode}
\end{figure}
\begin{lstlisting}[language=Python, caption=Python implementation]
from collections.abc import Iterator
def swap_elements(a_list,i,j):
	tmp = a_list[i];
	a_list[i] = a_list[j];
	a_list[j] = tmp    
class setperm(Iterator):
	def __init__(self,multiplicity):
		self.m = multiplicity
		self.k = len(multiplicity)
		self.P = []
		for i in range(self.k):
			self.P += [i+1]*multiplicity[i]
		self.n = len(self.P) 
		self.D = [1]*self.n
		self.T = 0 # No active element type           
	def __next__(self):
		if self.T == 0:
			self.T = 1
			return self.P.copy()
		else:
			return self.one_step(self.n).copy()
	def swap(self,i,j,df):
		swap_elements(self.P,i,j)
		swap_elements(self.D,i,j)
		for k in range(i+df,j,df):
			self.D[k] = 1
		self.T = 1
	def one_step(self,n): # One iteration of the 
		d = -1            # algorithm
		T = self.T
		for i in range(n-1,-1,-1): #
			if self.P[i] == T:
				d = i
				break
		df = self.D[d]
		j = d + df
		if d>-1:
			while j>-1 and j<self.n:
				if self.P[j] != T or self.D[j] != df:
					if self.P[j] > T:
						self.swap(d,j,df)
						return self.P
					break
				j = j+df
			self.D[d] = -self.D[d]
			return self.one_step(d)
		else: # No elements of type T can move!
			self.T = self.T+1 # Next type!
			if self.T >= self.k: # No elements can move.
				raise StopIteration # Exit!
			return self.one_step(self.n)
# Example code:
for perm in setperm([2,2,1,1]):
	print( perm )
\end{lstlisting}
\bibliographystyle{plainurl}


\begin{thebibliography}{10}

\bibitem{agaoka85}
Y.~Agaoka.
\newblock On the curvature of {Riemannian} submanifolds of codimension 2.
\newblock {\em Hokkaido Mathematical Journal}, 14:107--135, 1985.

\bibitem{Eades:1984}
Peter Eades and Brendan McKay.
\newblock An algorithm for generating subsets of fixed size with a strong
  minimal change property.
\newblock {\em Information Processing Letters}, 19(3):131--133, 1984.
\newblock URL:
  \url{https://www.sciencedirect.com/science/article/pii/0020019084900917},
  \href {https://doi.org/https://doi.org/10.1016/0020-0190(84)90091-7}
  {\path{doi:https://doi.org/10.1016/0020-0190(84)90091-7}}.

\bibitem{Selmer_M_Johnson:1963}
Selmer~M. Johnson.
\newblock Generation of permutations by adjacent transposition.
\newblock {\em Mathematics of Computation}, 17(83):282--285, 1963.
\newblock URL: \url{http://www.jstor.org/stable/2003846}.

\bibitem{Knuth:2005:ACPb}
Donald~E. Knuth.
\newblock {\em The Art of Computer Programming, Volume 4, Fascicle 3:
  Generating All Combinations and Partitions}.
\newblock Addison-Wesley Professional, 2005.

\bibitem{Korsh:2000}
James~F. Korsh and Paul LaFollette.
\newblock Multiset permutations and loopless generation of ordered trees with
  specified degree sequence.
\newblock {\em Journal of Algorithms}, 34(2):309--336, 2000.
\newblock URL:
  \url{https://www.sciencedirect.com/science/article/pii/S0196677499910593},
  \href {https://doi.org/https://doi.org/10.1006/jagm.1999.1059}
  {\path{doi:https://doi.org/10.1006/jagm.1999.1059}}.

\bibitem{Kruchinin:2022}
Vladimir Kruchinin, Yuriy Shablya, Dmitry Kruchinin, and Victor Rulevskiy.
\newblock Unranking small combinations of a large set in co-lexicographic
  order.
\newblock {\em Algorithms}, 15(2), 2022.

\bibitem{Rota1969}
Jacob T.~Schwartz Mark~Kac, Gian-Carlo~Rota.
\newblock {\em Discrete Thoughts, Essays on Mathematics, Science and
  Philosophy}.
\newblock Birkhäuser Boston, MA, 2 edition, 1992.

\bibitem{Mutze:2022}
Torsten Mütze.
\newblock Combinatorial gray codes-an updated survey, 2022.
\newblock URL: \url{https://arxiv.org/abs/2202.01280}, \href
  {https://doi.org/10.48550/ARXIV.2202.01280}
  {\path{doi:10.48550/ARXIV.2202.01280}}.

\bibitem{Ruskey2003}
Frank Ruskey.
\newblock Combinatorial generation.
\newblock Preliminary working draft. University of Victoria, Victoria, BC,
  Canada, Oct 2003.

\bibitem{Ruskey2005}
Frank Ruskey and Aaron Williams.
\newblock Generating combinations by prefix shifts.
\newblock volume 3595, pages 570--576, 08 2005.
\newblock \href {https://doi.org/10.1007/11533719_58}
  {\path{doi:10.1007/11533719_58}}.

\bibitem{Ruskey2009}
Frank Ruskey and Aaron Williams.
\newblock The coolest way to generate combinations.
\newblock {\em Discrete Mathematics}, 309(17):5305--5320, 2009.
\newblock Generalisations of de Bruijn Cycles and Gray Codes/Graph
  Asymmetries/Hamiltonicity Problem for Vertex-Transitive (Cayley) Graphs.
\newblock URL:
  \url{https://www.sciencedirect.com/science/article/pii/S0012365X07009570},
  \href {https://doi.org/https://doi.org/10.1016/j.disc.2007.11.048}
  {\path{doi:https://doi.org/10.1016/j.disc.2007.11.048}}.

\bibitem{Sawada2021}
Joe Sawada and Aaron Williams.
\newblock A universal cycle for strings with fixed-content (which are also
  known as multiset permutations).
\newblock 04 2021.

\bibitem{Williams2022}
Joe Sawada and Aaron Williams.
\newblock Constructing the first (and coolest) fixed-content universal cycle.
\newblock {\em Algorithmica}, pages 1--32, 11 2022.
\newblock \href {https://doi.org/10.1007/s00453-022-01047-2}
  {\path{doi:10.1007/s00453-022-01047-2}}.

\bibitem{Takaoka1999_a}
Tadao Takaoka.
\newblock An o(1) time algorithm for generating multiset permutations.
\newblock In {\em Algorithms and Computation}, pages 237--246, Berlin,
  Heidelberg, 1999. Springer Berlin Heidelberg.

\bibitem{Takaoka1999}
Tadao Takaoka.
\newblock {O(1) Time Algorithms for Combinatorial Generation by Tree
  Traversal}.
\newblock {\em The Computer Journal}, 42(5):400--408, 01 1999.

\bibitem{Takaoka2006}
Tadao Takaoka and Stephen Violich.
\newblock Combinatorial generation by fusing loopless algorithms.
\newblock In {\em Proceedings of the 12th Computing: The Australasian Theroy
  Symposium - Volume 51}, CATS '06, page 69–77, AUS, 2006. Australian
  Computer Society, Inc.

\bibitem{Torres2019}
Jose Torres-Jimenez and Idelfonso Izquierdo-Marquez.
\newblock A low spatial complexity algorithm to generate combinations with the
  strong minimal change property.
\newblock {\em Discrete Mathematics, Algorithms and Applications},
  11(05):1950060, 2019.
\newblock \href {https://doi.org/10.1142/S1793830919500605}
  {\path{doi:10.1142/S1793830919500605}}.

\bibitem{Trotter:1962}
H.~F. Trotter.
\newblock Algorithm 115: Perm.
\newblock {\em Commun. ACM}, 5(8):434–435, aug 1962.
\newblock \href {https://doi.org/10.1145/368637.368660}
  {\path{doi:10.1145/368637.368660}}.

\bibitem{Williams2009}
Aaron Williams.
\newblock Loopless generation of multiset permutations using a constant number
  of variables by prefix shifts.
\newblock pages 987--996, 01 2009.
\newblock \href {https://doi.org/10.1145/1496770.1496877}
  {\path{doi:10.1145/1496770.1496877}}.

\end{thebibliography}

\end{document}